\documentclass[11pt]{article}

\usepackage{latexsym}
\usepackage[dvips]{graphicx}
\usepackage{float}
\usepackage{amsmath}
\usepackage{epic}
\usepackage{color}
\usepackage{setspace}
\usepackage{lscape}
\usepackage{arydshln}

\begin{document}
\addtocounter{MaxMatrixCols}{20}
\bibliographystyle{plain}
\floatplacement{table}{H}
\newtheorem{definition}{Definition}[section]
\newtheorem{lemma}{Lemma}[section]
\newtheorem{theorem}{Theorem}[section]
\newtheorem{corollary}{Corollary}[section]
\newtheorem{proposition}{Proposition}[section]
\newcommand{\sni}{\sum_{i=1}^{n}}
\newcommand{\snj}{\sum_{j=1}^{n}}
\newcommand{\smj}{\sum_{j=1}^{m}}
\newcommand{\sumjm}{\sum_{j=1}^{m}}
\newcommand{\bdis}{\begin{displaymath}}
\newcommand{\edis}{\end{displaymath}}
\newcommand{\beq}{\begin{equation}}
\newcommand{\eeq}{\end{equation}}
\newcommand{\beqn}{\begin{eqnarray}}
\newcommand{\eeqn}{\end{eqnarray}}
\newcommand{\qed}{{\large $\sqcap$ \hskip -0.37cm $\sqcup$}}
\newcommand{\defeq}{\stackrel{\triangle}{=}}
\newcommand{\sep}{\;\;\;\;\;\; ; \;\;\;\;\;\;}
\newcommand{\real}{\mbox{$ I \hskip -4.0pt R $}}
\newcommand{\complex}{\mbox{$ I \hskip -6.8pt C $}}
\newcommand{\integ}{\mbox{$ Z $}}
\newcommand{\realn}{\real ^{n}}
\newcommand{\sqrn}{\sqrt{n}}
\newcommand{\sqrtwo}{\sqrt{2}}
\newcommand{\prf}{{\bf Proof. }}
\newcommand{\onehlf}{\frac{1}{2}}
\newcommand{\thrhlf}{\frac{3}{2}}
\newcommand{\fivhlf}{\frac{5}{2}}
\newcommand{\onethd}{\frac{1}{3}}
\newcommand{\lb}{\left ( }
\newcommand{\lcb}{\left \{ }
\newcommand{\lsb}{\left [ }
\newcommand{\labs}{\left | }
\newcommand{\rb}{\right ) }
\newcommand{\rcb}{\right \} }
\newcommand{\rsb}{\right ] }
\newcommand{\rabs}{\right | }
\newcommand{\lnm}{\left \| }
\newcommand{\rnm}{\right \| }
\newcommand{\lambdab}{\bar{\lambda}}
%
%
\newcommand{\xj}{x_{j}}
\newcommand{\xjb}{\bar{x}_{j}}
\newcommand{\xro}{x_{\rho}}
\newcommand{\xrob}{\bar{x}_{\rho}}
\newcommand{\xsig}{x_{\mu}}
\newcommand{\xsigb}{\bar{x}_{\mu}}
\newcommand{\xnmjb}{\bar{x}_{n-j+1}}
\newcommand{\xnmj}{x_{n-j+1}}
\newcommand{\aroj}{a_{\rho j}}
\newcommand{\arojb}{\bar{a}_{\rho j}}
\newcommand{\aroro}{a_{\rho \rho}}
\newcommand{\amuro}{a_{\mu \rho}}
\newcommand{\amumu}{a_{\mu \mu}}
\newcommand{\aii}{a_{ii}}
\newcommand{\aik}{a_{ik}}
\newcommand{\akj}{a_{kj}}
\newcommand{\atwoii}{a^{(2)}_{ii}}
\newcommand{\atwoij}{a^{(2)}_{ij}}
\newcommand{\ajj}{a_{jj}}
\newcommand{\aiib}{\bar{a}_{ii}}
\newcommand{\ajjb}{\bar{a}_{jj}}
\newcommand{\bii}{a_{jj}}
\newcommand{\biib}{\bar{a}_{jj}}
\newcommand{\aij}{a_{i,n-i+1}}
\newcommand{\akl}{a_{j,n-j+1}}
\newcommand{\aijb}{\bar{a}_{i,n-i+1}}
\newcommand{\aklb}{\bar{a}_{j,n-j+1}}
\newcommand{\bij}{a_{n-j+1,j}}
\newcommand{\arorob}{\bar{a}_{\rho \rho}}
\newcommand{\arosig}{a_{\rho \mu}}
\newcommand{\arosigb}{\bar{a}_{\rho \mu}}
\newcommand{\sumjrosig}{\sum_{\stackrel{j=1}{j\neq\rho,\mu}}^{n}}
\newcommand{\summuro}{\sum_{\stackrel{j=1}{j\neq\mu,\rho}}^{n}}
\newcommand{\sumjnoti}{\sum_{\stackrel{j=1}{j\neq i}}^{n}}
\newcommand{\sumlnoti}{\sum_{\stackrel{\ell=1}{\ell \neq i}}^{n}}
\newcommand{\sumknoti}{\sum_{\stackrel{k=1}{k\neq i}}^{n}}
\newcommand{\sumknotij}{\sum_{\stackrel{k=1}{k\neq i,j}}^{n}}
\newcommand{\sumk}{\sum_{k=1}^{n}}
\newcommand{\snl}{\sum_{\ell=1}^{n}}
\newcommand{\sumji}{\sum_{\stackrel{j=1}{j\neq i, n-i+1}}^{n}}
\newcommand{\sumki}{\sum_{\stackrel{k=1}{k\neq i, n-i+1}}^{n}}
\newcommand{\sumkj}{\sum_{\stackrel{k=1}{k\neq j, n-j+1}}^{n}}
\newcommand{\sumjro}{\sum_{\stackrel{j=1}{j\neq\rho}}^{n}}
\newcommand{\rrosig}{R''_{\rho \mu}}
\newcommand{\rro}{R'_{\rho}}
\newcommand{\gamror}{\Gamma_{\rho}^{R}(A)}
\newcommand{\gamir}{\Gamma_{i}^{R}(A)}
\newcommand{\gamctrr}{\Gamma_{\frac{n+1}{2}}^{R}(A)}
\newcommand{\gamctrc}{\Gamma_{\frac{n+1}{2}}^{C}(A)}
\newcommand{\gamroc}{\Gamma_{\rho}^{C}(A)}
\newcommand{\gamjc}{\Gamma_{j}^{C}(A)}
\newcommand{\lamror}{\Lambda_{\rho}^{R}(A)}
\newcommand{\lamir}{\Lambda_{i}^{R}(A)}
\newcommand{\lamirepsilon}{\Lambda_{i}^{R}(A_{\epsilon})}
\newcommand{\lamnir}{\Lambda_{n-i+1}^{R}(A)}
\newcommand{\lamjr}{\Lambda_{j}^{R}(A)}
\newcommand{\phiij}{\Phi_{ij}^{R}(A)}
\newcommand{\delir}{\Delta_{i}^{R}(A)}
\newcommand{\vir}{V_{i}^{R}(A)}
\newcommand{\pamir}{\Pi_{i}^{R}(A)}
\newcommand{\xir}{\Xi_{i}^{R}(A)}
\newcommand{\lamjc}{\Lambda_{j}^{C}(A)}
\newcommand{\vjc}{V_{j}^{C}(A)}
\newcommand{\pamjc}{\Pi_{j}^{C}(A)}
\newcommand{\xjc}{\Xi_{j}^{C}(A)}
\newcommand{\lamroc}{\Lambda_{\rho}^{C}(A)}
\newcommand{\lamsigr}{\Lambda_{\mu}^{R}(A)}
\newcommand{\lamsigc}{\Lambda_{\mu}^{C}(A)}
\newcommand{\psii}{\Psi_{i}^{R}(A)}
\newcommand{\psiq}{\Psi_{q}}
\newcommand{\psiiepsilon}{\Psi_{i}^{R}(A_{\epsilon})}
\newcommand{\psiqepsilon}{\Psi_{q}(A_{\epsilon})}
\newcommand{\psiqc}{\Psi_{q}^{c}}
\newcommand{\psiqcepsilon}{\Psi_{q}^{c}(A_{\epsilon})}

\newcommand{\xmu}{x_{\mu}}
\newcommand{\xmub}{\bar{x}_{\mu}}
\newcommand{\xnu}{x_{\nu}}
\newcommand{\xnub}{\bar{x}_{\nu}}
\newcommand{\amuj}{a_{\mu j}}
\newcommand{\amujb}{\bar{a}_{\mu j}}
\newcommand{\amumub}{\bar{a}_{\mu \mu}}
\newcommand{\amunu}{a_{\mu \nu}}
\newcommand{\amunub}{\bar{a}_{\mu \nu}}
\newcommand{\sumjmunu}{\sum_{\stackrel{j=1}{j\neq\mu,\nu}}}
\newcommand{\rmunu}{R''_{\mu \nu}}
\newcommand{\rmu}{R'_{\mu}}

\newcommand{\azero}{a_{0}}
\newcommand{\aone}{a_{1}}
\newcommand{\atwo}{a_{2}}
\newcommand{\ath}{a_{3}}
\newcommand{\afr}{a_{4}}
\newcommand{\afv}{a_{5}}
\newcommand{\asx}{a_{6}}
\newcommand{\anmo}{a_{n-1}}
\newcommand{\akmo}{a_{k-1}}
\newcommand{\akpo}{a_{k+1}}
\newcommand{\anmt}{a_{n-2}}
\newcommand{\anmth}{a_{n-3}}
\newcommand{\akmt}{a_{k-2}}
\newcommand{\an}{a_{n}}
\newcommand{\ak}{a_{k}}

\newcommand{\bzero}{b_{0}}
\newcommand{\bone}{b_{1}}
\newcommand{\btwo}{b_{2}}
\newcommand{\bth}{b_{3}}
\newcommand{\bfr}{b_{4}}
\newcommand{\bfv}{b_{5}}
\newcommand{\bsx}{b_{6}}
\newcommand{\bnmo}{b_{n-1}}
\newcommand{\bkmo}{b_{k-1}}
\newcommand{\bkpo}{b_{k+1}}
\newcommand{\bnmt}{b_{n-2}}
\newcommand{\bkmt}{b_{k-2}}
\newcommand{\bn}{b_{n}}
\newcommand{\bk}{b_{k}}

\newcommand{\betk}{\beta^{-k}}
\newcommand{\betkmo}{\beta^{1-k}}
\newcommand{\bo}{\beta^{-1}}
\newcommand{\betnmk}{\beta^{n-k}}

\newcommand{\ix}{x}

\newcommand{\alfazero}{\alpha_{0}}
\newcommand{\alfaone}{\alpha_{1}}
\newcommand{\alfatwo}{\alpha_{2}}
\newcommand{\alfath}{\alpha_{3}}
\newcommand{\alfafr}{\alpha_{4}}
\newcommand{\alfafv}{\alpha_{5}}
\newcommand{\alfasx}{\alpha_{6}}
\newcommand{\alfanmo}{\alpha_{n-1}}
\newcommand{\alfanmt}{\alpha_{n-2}}
\newcommand{\alfan}{\alpha_{n}}
\newcommand{\alfak}{\alpha_{k}}
\newcommand{\bkl}{b_{k\ell}}

\newcommand{\abii}{(AB)_{ii}}
\newcommand{\abil}{(AB)_{i\ell}}
\newcommand{\matrixspace}{\;\;}
\newcommand{\ellone}{\ell_{1}}
\newcommand{\elltwo}{\ell_{2}}

\newcommand{\resh}{\rho}
\newcommand{\snk}{s}
\newcommand{\ronenk}{rone}
\newcommand{\rtwonk}{rtwo}
\newcommand{\obar}{\bar{O}}
\newcommand{\Lbar}{\bar{L}}
\newcommand{\dltaone}{\delta_{1}}
\newcommand{\rhoone}{\rho_{1}}
\newcommand{\rhotwo}{\rho_{2}}
\newcommand{\rone}{r}
\newcommand{\rtwo}{R}

\newcommand{\ftilde}{\widetilde{f_{k}}}
\newcommand{\scr}{\mathcal{H}}
\newcommand{\gvulscr}{\partial \mathcal{R}}

\begin{center}
\large
{\bf GEOMETRIC ASPECTS OF PELLET'S AND RELATED THEOREMS} \\
\vskip 0.5cm
\normalsize
A. Melman \\
Department of Applied Mathematics \\
School of Engineering, Santa Clara University  \\
Santa Clara, CA 95053  \\
e-mail : amelman@scu.edu \\
\vskip 0.5cm
\end{center}

\begin{abstract}
Pellet's theorem determines when the zeros of a polynomial can be separated into two regions,
according to their moduli.
We refine one of those regions and replace it with the closed interior of a lemniscate that provides
more precise information on the location of the zeros.
Moreover, Pellet's theorem is considered
the generalization of a zero inclusion region due to Cauchy. Using linear algebra tools, we derive a different 
generalization that leads to a sequence of smaller inclusion regions, which are also the closed interiors of lemniscates. 
\vskip 0.15cm
{\bf AMS(MOS) subject classification :} 12D10, 15A18, 30C15
\end{abstract}


\section{Introduction}


Pellet's classical theorem (\cite{Pellet},\cite[Theorem (28,1)]{Marden})
derives a criterion for the separation of the zeros of a general polynomial with complex coefficients 
into two regions of the
complex plane: a disk and the complement of a larger disk, both centered at the origin. 
Our first result is to replace the latter set by the interior of a lemniscate, which provides more precise 
information on the location of the zeros. It adds a geometric
component to a theorem that is formulated in terms of simple bounds on the moduli of the zeros that obscure
the details of the zero distribution. Although Pellet's theorem is often viewed as the generalization of an
inclusion region by Cauchy (\cite{Cauchy},\cite[Theorem (27,1)]{Marden}), 
we show that a different generalization can be obtained, leading to smaller regions 
consisting, once again, of the closed interiors of lemniscates.

Both Pellet's theorem and Cauchy's result can be proven with Rouch\'{e}'s theorem. 
However, to refine and extend these theorems, we used 
the Gershgorin set to estimate the eigenvalues of a polynomials's companion matrix. 
This set is a union of disks in the complex plane, centered at the diagonal elements of the matrix. 
To be able to extract useful results from this often crude method, we consider a similarity transformation
of an appropriate polynomial of the companion matrix, rather than the companion matrix itself. 
The main advantage of the Gershgorin set, as we apply it here, is to point out results that might 
otherwise not be apparent, even if their subsequent proof by Rouch\'{e}'s theorem is relatively straightforward.

Although Pellet's theorem was recently generalized to matrix polynomials in~\cite{BiniMatPol} and~\cite{MelmanMatPol},
our results here cannot easily be similarly extended because of their heavy dependence on the scalar nature of
the coefficients. We will discuss other improvements after the derivation of our results.

A good introduction to Gershgorin disks and other eigenvalue inclusion regions can be found in \cite[Ch. 6]{HJ}. 
For a more in-depth study of the subject, including its interesting history, we refer to \cite{V} and the many 
references therein. For results concerning polynomial zeros we refer to the encyclopedic work \cite{Marden}.

In Section 2, we state the aforementioned theorems, together with definitions and lemmas that are needed in 
Section 3, where we derive and illustrate our main results.

\section{Preliminaries}


We start by stating Pellet's theorem.

\begin{theorem}
(\cite{Pellet}, \cite[Theorem (28,1), p.128]{Marden})
\label{Pellet}
Given the polynomial $p(z)=z^{n} + a_{n-1} z^{n-1} + \dots + \aone z + \azero$ with complex coefficients, 
$\azero \ak \neq 0$, and $n \geq 3$. Let $1 \leq k \leq n-1$, and let the polynomial 
\bdis
f_{k}(x)=x^{n} + |a_{n-1}| x^{n-1} + \dots + |a_{k+1}| x^{k+1} - |a_{k}| x^{k} + |a_{k-1}|x^{k-1} + \dots  + |\azero| 
\edis
have two distinct positive roots $r$ and $R$, $r < R$.
Then $p$ has exactly $k$ zeros in or on the circle $|z| = r$ and no zeros in the annular ring $r < |z| < R$.
\end{theorem}

We will not dwell on the numerical solution of $f_{k}(x)=0$.
A systematic method to do so can be found in~\cite{MelmanNumAlg}, while a heuristic method was developed 
in~\cite{Rump}. An implementation using the Newton polygon can be found in~\cite{BiniNumAlg}
and~\cite{BiniMatPol} for scalar and matrix polynomials, respectively.

Pellet's theorem is considered the generalization of the following result by Cauchy:

\begin{theorem}
(\cite{Cauchy}, \cite[Theorem (27,1), p.122]{Marden})
\label{Cauchy}
All the zeros of the polynomial $p(z)=z^{n} + a_{n-1} z^{n-1} + \dots + \aone z + \azero$ 
with complex coefficients lie in the circle $|z| = r$, where $r$ is the positive root of the equation
\bdis
x^{n} - |a_{n-1}| x^{n-1} - \dots - |\aone| x - |\azero| = 0 \; . 
\edis
\end{theorem}

Both these theorems are a direct consequence of Rouch\'{e}'s theorem (\cite{Rouche}, \cite[Theorem 1.6, p.181]{Lang}). 
However, they only provide information on the moduli of the zeros. To introduce more interesting geometry into these 
results, leading to better information on the distribution of the zeros in the complex plane,
we will instead use linear algebra tools, namely, Gershgorin's theorem and the polynomial's companion 
matrix. We state Gershgorin's theorem next.

%
%
\begin{theorem}
\label{Gershgorin}
(Gershgorin, \cite{G}, \cite[Theorem 6.1.1, p.344]{HJ})
All the eigenvalues of the $n \times n$ complex matrix $A$ with elements $a_{ij}$ are located in the union of
$n$ disks
\bdis
\Gamma(A) = \bigcup_{i=1}^{n} \lcb z \in \complex \, : \, |z - a_{ii} | \leq R'_{i}(A) \rcb , \; \text{with} \; 
R'_{i}(A) = \sum_{\stackrel{j=1}{j\neq i}}^{n} |a_{ij} | \; .
\edis
Moreover, if $\ell$ disks form a connected region that is disjoint from the remaining $n-\ell$ disks, then this region 
contains exactly $\ell$ eigenvalues.
\end{theorem}

$R'_{i}(A)$ is called the $i$th deleted row sum of $A$. 
The spectrum of $A$ and $A^{T}$ is the same, so that the Gershgorin set also has a column version, 
obtained by applying the theorem to $A^{T}$, where the deleted column sums replace the deleted row sums. 
In addition, any similarity transformation of $A$, namely, $S^{-1}AS$ for a nonsingular matrix $S$, has the
same eigenvalues as $A$, but may have a smaller Gershgorin set. Frequently, $S$ is chosen to be a diagonal matrix. 

Eigenvalue inclusion sets can be used to estimate zeros of a polynomial by applying them to the 
polynomial's companion matrix, whose eigenvalues are the zeros of the polynomial. 
A common choice for a companion matrix of the monic polynomial 
$p(z)=z^{n} + \anmo z^{n-1} + \dots + \aone z + \azero$ 
is given by (see, e.g., \cite[p.146]{HJ}):
\beq
\label{companionmatrix}
C(p) =
\begin{pmatrix}
0 &        &       &   & -\azero     \\
1 &        &       &   & -\aone     \\
  & \ddots &       &   & \vdots \\
  &        &       & 1 & -\anmo    \\
\end{pmatrix}
\; ,
\eeq
where blank entries represent zeros, a convention we will follow throughout.
In what follows, we set $a_{n}=1$, where $a_{n}$ is the leading coefficient of the aforementioned polynomial $p$,
and denote an open disk centered at $a$ with radius $\rho$ by $O(a;\rho)$.
The closure and the complement of a set $\Delta$ will be denoted by $\bar{\Delta}$ and $\Delta^{c}$,
respectively. We also define the following.

\begin{definition}
The associated polynomials $\{ p_{k} \}_{k=1}^{n-1}$ of the polynomial 
\bdis
p(z)=z^{n} + a_{n-1} z^{n-1} + \dots + a_{0}
\edis
are defined by $p_{1}(z)=z$ and the recursion
\bdis
p_{k+1}(z) = z \lb p_{k}(z) + a_{n-k} \rb \; .
\edis
\end{definition}
Therefore, $p_{k}$, $1 \leq k \leq n-1$, is given by
\bdis
p_{k}(z)=z^{k} + a_{n-1} z^{k-1} + \dots + a_{n-k+1}z \; .
\edis
\begin{definition}
The polynomial $P_{k}$ 
is obtained from $p_{k}$ by replacing its coefficients with their moduli, i.e.,
\bdis
P_{k}(z)=z^{k} + |a_{n-1}| z^{k-1} + \dots + |a_{n-k+1}|z \; .
\edis
\end{definition}
\begin{definition}
The complex $n \times n$ matrix $M_{k}(p)$, $1 \leq k \leq n-1$, is defined as 
\bdis
M_{k}(p) = p_{k}(C(p)) \; ,
\edis
where $C(p)$ is the companion matrix of $p$.
\end{definition}
The eigenvalues of $M_{k}(p)$ are $\{ p_{k}(z_{i}) \}_{i=1}^{n}$, where
$\{ z_{i} \}_{i=1}^{n}$ are the eigenvalues of $C(p)$, which are also the zeros of $p$. 
Its structure is derived in the following lemma.

\begin{lemma}
\label{lemmaMk}
The matrix $M_{k}(p)$, $1 \leq k \leq n-1$, is given by      
\bdis
M_{k}(p) = 
\begin{pmatrix}
 0     &         &        &         & -\azero      &               &             &           &           \\
a_{n-k+1}&  \ddots &        &         & -\aone       & \ddots        &             &           &           \\
a_{n-k+2}&  \ddots & \ddots &         &  \vdots      & \ddots        & \ddots      &           &           \\
\vdots &  \ddots & \ddots &    0    & -a_{n-k-1}     & \ddots        & \ddots      & -\azero   &           \\
a_{n-1}& \ddots  & \ddots & a_{n-k+1} & -a_{n-k}       & \ddots        & \ddots      & -\aone    &  -\azero  \\
1      & \ddots  & \ddots & a_{n-k+2} &              & \ddots        & \ddots      &  \vdots   &  -\aone   \\
       & \ddots  & \ddots & \vdots  &              &               & \ddots      & -a_{n-k-1}  &   \vdots  \\
       &         & \ddots & a_{n-1} &              &               &             & -a_{n-k}    & -a_{n-k-1}  \\
       &         &        & 1       &              &               &             &           & -a_{n-k}    \\
\end{pmatrix}
\; ,           
\edis
where $k$ diagonal elements are equal to $a_{n-k}$, while the remaining ones are zero.
\end{lemma}
\prf The proof is by induction. Since $M_{1}(p) = C(p)$, the lemma is obviously true for $k=1$. Now assume that is true
for $M_{j}(p)$, $1 \leq j \leq n-2$. A straightforward calculation then shows that $C(p) M_{j}(p)$
is given by
\bdis
\begin{pmatrix}
0         &         &        &         &              &                & -\azero      &               &               &             &             &  a_{n-j}\azero               \\
0         &  \ddots &        &         &              &                & -\aone       &    \ddots     &               &             &             &  a_{n-j}\aone                \\
a_{n-j+1} &  \ddots & \ddots &         &              &                & -\atwo       &    \ddots     & \ddots        &             &             &  a_{n-j}\atwo                    \\
a_{n-j+2} &  \ddots & \ddots & \ddots  &              &                & \vdots       &    \ddots     & \ddots        & \ddots      &             &  \vdots                    \\
\vdots    & \ddots  & \ddots & \ddots  &  \ddots      &                & \vdots       &    \ddots     & \ddots        & \ddots      & -\azero     & a_{n-j}a_{n-j-3}               \\
a_{n-1}   & \ddots  & \ddots & \ddots  &  \ddots      &  0             & -a_{n-j-2}   &    \ddots     & \ddots        & \ddots      & -\aone      & a_{n-j}a_{n-j-2} -\azero       \\
1         & \ddots  & \ddots & \ddots  &  \ddots      &  0             & -a_{n-j-1}   &    \ddots     & \ddots        & \ddots      & -\atwo      & a_{n-j}a_{n-j-1}  -\aone                   \\
          & \ddots  & \ddots & \ddots  &  \ddots      &  a_{n-j+1}     & -a_{n-j}     &    \ddots     & \ddots        & \ddots      & \vdots      &   \vdots                   \\
          &         & \ddots & \ddots  &  \ddots      &  a_{n-j+2}     &              &    \ddots     & \ddots        & \ddots      & \vdots      &  \vdots                        \\
          &         &        & \ddots  &  \ddots      &  \vdots        &              &               & \ddots        & \ddots      & -a_{n-j-2}  & a_{n-j}a_{n-3}   - a_{n-j-3}          \\
          &         &        &         &  \ddots      &  a_{n-1}       &              &               &               & \ddots      & -a_{n-j-1}  & a_{n-j}a_{n-2}   - a_{n-j-2}   \\
          &         &        &         &              &  1             &              &               &               &             & -a_{n-j}    & a_{n-j}a_{n-1}   - a_{n-j-1}   \\
\end{pmatrix}
\; ,
\edis 
from which one easily deduces that $M_{j+1}(p)=C(p) M_{j}(p) + a_{n-j}C(p)$ is of the same form as $M_{j}(p)$. \qed

The form of $M_{k}(p)$ makes it convenient to apply the column version of Gershgorin's theorem since the deleted column sums
are the same for identical diagonal elements. To add flexibility to the Gershgorin set, we will use a diagonal
similarity transformation. The next lemma shows its effect on $M_{k}(p)$.

\begin{lemma}
\label{lemmaMksimilarity}
Let $D_{\ix}$ be a diagonal matrix with diagonal $\lb \ix^{n}, \ix^{n-1}, \dots, \ix \rb$ 
for $\ix > 0$.
Then the matrix $D^{-1}_{\ix} M_{k}(p) D_{\ix}$, $1 \leq k \leq n-1$, is given by 
\bdis
\begin{pmatrix}
 0                    &         &        &                       & -\azero /\ix^{n-k}        &               &             &                         &                        \\
a_{n-k+1} \ix         &  \ddots &        &                       & -\aone / \ix^{n-k-1}      & \ddots        &             &                         &                        \\
a_{n-k+2} \ix^{2}     &  \ddots & \ddots &                       & \vdots                    & \ddots        & \ddots      &                         &                        \\
\vdots                &  \ddots & \ddots &    0                  & -a_{n-k-1} / \ix          & \ddots        & \ddots      & -\azero/\ix^{n-k}       &                        \\
a_{n-1} \ix^{k-1} & \ddots  & \ddots & a_{n-k+1} \ix         & -a_{n-k}                    & \ddots        & \ddots      & -\aone/ \ix^{n-k-1}     &  -\azero/\ix^{n-k}     \\
\ix^{k}           & \ddots  & \ddots & a_{n-k+2} \ix^{2}     &                           & \ddots        & \ddots      & \vdots                  &  -\aone/ \ix^{n-k-1}   \\
                      & \ddots  & \ddots & \vdots                &                           &               & \ddots      & -a_{n-k-1}/ \ix         &  \vdots                \\
                      &         & \ddots & a_{n-1} \ix^{k-1} &                           &               &             & -a_{n-k}                  & -a_{n-k-1}/ \ix        \\
                      &         &        & \ix^{k}           &                           &               &             &                         & -a_{n-k}                 \\
\end{pmatrix}
\; \cdot
\edis
\end{lemma}
\prf For any diagonal matrix $D$ with diagonal $(d_{1},d_{2},...,d_{n})$, and matrix $A$ with elements $a_{ij}$,
$(D^{-1}AD)_{ij} = d_{j}a_{ij}/d_{i}$. The lemma then follows directly
by substituting $D_{\ix}$ in $D^{-1}_{\ix} M_{k}(p) D_{\ix}$.~\qed

In what follows we will frequently encounter lemniscates of the form $|q(z)|=\alpha$, where $q$ is a polynomial.
The zeros of $q$ are the foci of the lemniscate, which, depending on the value of $\alpha$ can consist of at most
$m$ disjoint closed curves, where $m$ is the order of $q$. These curves are simple except for at most 
$m-1$ critical values of $\alpha$. 
If a lemniscate has distinct foci $\{z_{j}\}_{j=1}^{m}$, and if there exists $\eta > 0$ so that
the disks $\Delta_{j} = \lcb z \in \complex \, : \, |z - z_{j}| < \eta \rcb$ are disjoint, then the lemniscate
$ \lcb z \in \complex \, : \,  \Pi_{j=1}^{m} |z - z_{j}|  = \rho^{n} \rcb$
is contained in the union $\cup_{j=1}^{m} \Delta_{j}$ for any $\rho \leq \eta$.
Since a lemniscate must contain all of its foci in its interior, this provides an easily computable sufficient 
condition for a lemniscate to be composed of $m$ disjoint simple curves.
We refer to \cite[Vol. I, p.379]{Markushevich} for a more detailed discussion of lemniscates.


\section{Main results}


Our first result is Pellet's theorem with a refinement of the region outside the disk with radius $R$ in 
Theorem~\ref{Pellet}, where the largest zeros can be found.

\begin{theorem}
\label{geometricPellet}
Let $p(z)=z^{n} + a_{n-1} z^{n-1} + \dots + \aone z + \azero$ be a polynomial with complex coefficients, 
$\azero \neq 0$, $1 \leq k \leq n-1$, $n \geq 3$, and with zeros $\{z_{i}\}_{i=1}^{n}$, labeled so that 
$|z_{1}| \leq |z_{2}| \leq \dots \leq |z_{n}|$.
Let $\{p_{j}\}_{j=1}^{n-1}$ be the associated polynomials of $p$, 
let $P_{j}$ be the polynomial obtained from $p_{j}$ by replacing its coefficients with their moduli, and let
\bdis
\mu(k,x)=\sum_{j=0}^{k-1} |a_{j}|x^{j-k} \; .  
\edis
Furthermore, let
\bdis
f_{k}(x)=x^{n} + |a_{n-1}| x^{n-1} + \dots + |a_{k+1}| x^{k+1} - |a_{k}| x^{k} + |a_{k-1}|x^{k-1} + \dots  + |\azero| 
\edis
have two distinct positive roots $r$ and $R$ such that $0 < \rone < \rtwo$, and define the disjoint sets              
\begin{eqnarray*}
& & \Omega_{1}(k) = \lcb z \in \complex \, : \, |p_{n-k}(z) | \leq P_{n-k}(\rone) = |\ak| - \mu(k,r) \rcb  \; , \\    
& & \Omega_{2}(k) = \lcb z \in \complex \, : \, |p_{n-k}(z) + a_{k}| \leq \mu(k,\rtwo) = |\ak| - P_{n-k}(R) \rcb \; .
\end{eqnarray*}
Then:
\newline {\bf (1)} the $k$ zeros $\{z_{i}\}_{i=1}^{k}$ are contained in the closed disk $\obar(0;\rone)$, whereas 
the remaining $n-k$ zeros $\{z_{i}\}_{i=k+1}^{n}$ are contained in $\Omega_{2}(k)$, the closed interior of
a lemniscate;
\newline {\bf (2)} $\obar(0;\rone) \subseteq \Omega_{1}(k)$ and $\Omega_{2}(k) \subseteq O^{c}(0;\rtwo)$; 
\newline {\bf (3)} if $\Omega_{2}(k)$ consists of disjoint regions whose boundaries are simple closed (Jordan) curves, 
then each disjoint region contains as many
zeros of $p$ as it contains zeros of $p_{n-k}(z) + a_{k}$ (or foci of $\Omega_{2}(k)$).
\end{theorem}
\prf 
We start by observing that $\ak \neq 0$ since $f_{k}$ has positive zeros, and that $f_{k}$ can be written as
\beq
\label{fmueq}
f_{k}(x) = x^{k} \lb P_{n-k}(x) - |\ak| + \mu(k,x) \rb \; .
\eeq
That the zeros $\{z_{i}\}_{i=1}^{k}$ are contained in the closed disk $\obar(0;\rone)$ is obtained from 
Rouch\'{e}'s theorem, exactly as in the proof of Pellet's theorem (see, e.g., \cite[Theorem (28,1), p. 128]{Marden}). 

The numbers $\{p_{n-k}(z_{i}) \}_{i=1}^{n}$ are the eigenvalues of $p_{n-k}(C(p))=M_{n-k}(p)$, and therefore also
of $D^{-1}_{\ix} M_{n-k}(p) D_{\ix}$ for any $x > 0$, where $D_{\ix}$ is as in Lemma~\ref{lemmaMksimilarity}. 
From that same lemma, the matrix $D^{-1}_{\ix} M_{n-k}(p) D_{\ix}$ is given by
\bdis
\begin{pmatrix}
 0                    &         &        &                       & -\azero /\ix^{k}        &               &             &                         &                        \\
a_{k+1} \ix         &  \ddots &        &                       & -\aone / \ix^{k-1}      & \ddots        &             &                         &                        \\
a_{k+2} \ix^{2}     &  \ddots & \ddots &                       & \vdots                    & \ddots        & \ddots      &                         &                        \\
\vdots                &  \ddots & \ddots &    0                  & -a_{k-1} / \ix          & \ddots        & \ddots      & -\azero/\ix^{k}       &                        \\
a_{n-1} \ix^{n-k-1} & \ddots  & \ddots & a_{k+1} \ix         & -a_{k}                    & \ddots        & \ddots      & -\aone/ \ix^{k-1}     &  -\azero/\ix^{k}     \\
\ix^{n-k}           & \ddots  & \ddots & a_{k+2} \ix^{2}     &                           & \ddots        & \ddots      & \vdots                  &  -\aone/ \ix^{k-1}   \\
                      & \ddots  & \ddots & \vdots                &                           &               & \ddots      & -a_{k-1}/ \ix         &  \vdots                \\
                      &         & \ddots & a_{n-1} \ix^{n-k-1} &                           &               &             & -a_{k}                  & -a_{k-1}/ \ix        \\
                      &         &        & \ix^{n-k}           &                           &               &             &                         & -a_{k}                 \\
\end{pmatrix}
\; ,
\edis
and its Gershgorin column set is easily seen to be the union of a disk centered at the origin with radius 
$P_{n-k}(x)$ and a disk centered at $-a_{k}$ with radius $\mu(k,x)$. 
These disks are disjoint if there exists $\delta > 0$ such that
$P_{n-k}(\delta) + \mu(k,\delta) < |a_{k}|$, which, in view of~(\ref{fmueq}), is equivalent to $f_{k}(\delta) < 0$. 
It therefore suffices to choose any $\delta$ for which $\rone < \delta < \rtwo$ to obtain disjoint disks. 
Since the diagonal of $D^{-1}_{\delta} M_{n-k}(p) D_{\delta}$ contains $k$ zeros, 
we conclude that exactly $k$ of the $n$ numbers $\{p_{n-k}(z_{i}) \}_{i=1}^{n}$ lie in the closed disk 
$\obar \lb 0;P_{n-k}(\delta) \rb $, while the remaining $n-k$ lie in the disjoint closed disk
$\obar \lb -a_{k};\mu(k,\delta) \rb $.
Since this is true for any $\delta$ such that $\rone < \delta < \rtwo$, the same conclusions hold for
the disks
$\obar \lb 0;P_{n-k}(r) \rb $ and $\obar \lb -a_{k};\mu(k,\rtwo) \rb $, respectively.
This concludes the proof of the first part of the theorem.

Recalling that we defined $a_{n}=1$, we now observe that $|z| \leq \rone$ implies that 
\beq
\label{firstkineq}
\left | p_{n-k}(z) \right | 
=    \left | \sum_{j=k+1}^{n} a_{j} z^{j-k} \right | 
\leq \sum_{j=k+1}^{n} |a_{j}| |z|^{j-k}  
\leq \sum_{j=k+1}^{n} |a_{j}| \rone^{j-k} 
=    P_{n-k}(\rone) \; ,
\eeq
which means that $\obar (0;r) \subseteq \Omega_{1}(k)$.
It also means that it is the $k$ numbers $\{p_{n-k}(z_{i}) \}_{i=1}^{k}$, corresponding to the first $k$ zeros of $p$,
that lie in $\obar \lb 0;P_{n-k}(\rone) \rb $. 
To show the second inclusion in the statement of the theorem, assume that $z \in \Omega_{2}(k)$. Since the disk 
$\obar \lb -a_{k};\mu(k,\rtwo) \rb$ is bounded away from the origin, we have that  
$|p_{n-k}(z)| \geq |a_{k}|-\mu(k,\rtwo) = P_{n-k}(R) > 0$, and therefore that
\beq
\label{secondkineq}
P_{n-k}(|z|) \geq |p_{n-k}(z)| \geq P_{n-k}(R) > 0 \; .
\eeq
The polynomial $P_{n-k}$ is strictly increasing for positive arguments, so that inequality~(\ref{secondkineq})
implies that $|z| \geq R$ and therefore that $z$ must lie in $O^{c}(0;R)$.
This proves the second part of the theorem.

The lemniscate $\Omega_{2}(k)$ has $n-k$ foci, which are the zeros of $p_{n-k}(z) + a_{k}$, and it can consist of,
at most, $n-k$ disjoint regions. Any disjoint region contains one or more foci. 
Let us now assume that there exists a disjoint region $\scr$ of $\Omega_{2}(k)$ with a simple boundary. 
Define 
\bdis
\widetilde{\Omega}_{2}(k) = \lcb z \in \complex \, : \, |p_{n-k}(z) + a_{k}| \leq \mu(k,\delta) \rcb \; ,
\; r < \delta < R \; ,
\edis
where $R-\delta$ is small enough so that $\widetilde{\Omega}_{2}(k)$ (which contains
$\Omega_{2}$ since $\mu(k,\delta)~>~\mu(k,R)$) has the same number
of disjoint regions with simple boundary as $\Omega_{2}$, one of which, $\widetilde{{\scr}}$, contains ${\scr}$.
Similarly as before, we have for $z \in \widetilde{\Omega}_{2}(k)$ 
that $|p_{n-k}(z)| \geq |a_{k}|-\mu(k,\delta) > P_{n-k}(\delta) > 0$, and therefore that
\bdis
P_{n-k}(|z|) \geq |p_{n-k}(z)| > P_{n-k}(\delta) > 0 \; ,
\edis
so that now $|z| > \delta$.

For any $z$ on $\partial\widetilde{{\scr}}$, we have that
\beq
\label{lemdef}
\left | p_{n-k}(z) + a_{k} \right | = \mu(k,\delta) \; .
\eeq
Since $|z| > \delta$, multiplying both sides of~(\ref{lemdef}) by $|z|^{k}$ yields
\begin{eqnarray*}
\left | z^{n} + a_{n-1} z^{n-1} + \dots + a_{k}z^{k} \right | 
         &   = & |a_{k-1}| \left | \dfrac{z}{\delta} \right | |z|^{k-1}
                 + |a_{k-2}| \left | \dfrac{z}{\delta} \right |^{2} |z|^{k-2}
                  + \dots + |a_{0}| \left | \dfrac{z}{\delta} \right |^{k} \\ 
         &   >  & |a_{k-1}| |z|^{k-1} + |a_{k-2}| |z|^{k-2} + \dots + |a_{0}|  \\ 
         & \geq & \left | a_{k-1} z^{k-1} + a_{k-2} z^{k-2} + \dots + a_{0} \right | \; .
\end{eqnarray*}
Then, by Rouch\'{e}'s theorem, the polynomial $p$ has as many zeros inside this disjoint region $\widetilde{{\scr}}$ as 
$z^{k}(p_{n-k}(z)+a_{k})$, which is as many zeros as $p_{n-k}(z)+a_{k}$ since $\widetilde{{\scr}}$ does not contain 0.
Because this remains true as $\delta \rightarrow R^{-}$, the proof follows. \qed

It is a direct consequence of part (3) of this theorem that if $\Omega_{2}(k)$ consists of $n-k$ disjoint
regions with simple boundaries, then each region must contain exactly one zero of $p$.

\begin{figure}[H]
\begin{center}
\includegraphics[width=0.45\linewidth]{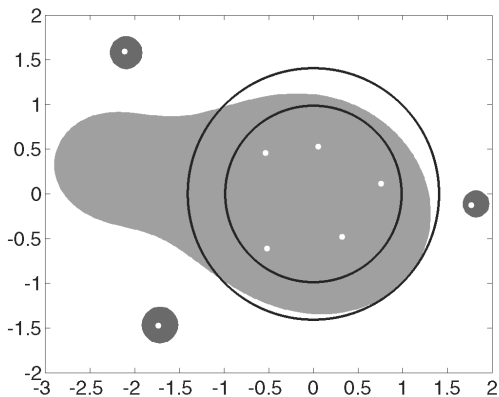}
\caption{The sets $\Omega_{1}(5)$ (light gray) and $\Omega_{2}(5)$ (dark gray) for $q_{A}$.}
\label{Pelletfig}   
\end{center}
\end{figure}

We illustrate this theorem at the hand of the polynomial $q_{A}$, defined by
\bdis
q_{A}(z) = z^{8} + 2 z^{7} - (1-i) z^{6} - 11 z^{5} + \dfrac{1}{2} z^{4} + z^{3} + \dfrac{5}{2} z^2 - z + (1+i) \; ,
\edis
and its sets $\Omega_{1}(5)$ and $\Omega_{2}(5)$. We obtain $r=0.9872$ and $R=1.4065$, so that
\begin{eqnarray*}
& & \Omega_{1}(5) = \lcb z \in \complex \, : \, |z^{3} + 2 z^{2} - (1-i) z| \leq 4.3065 \rcb  \; , \\    
& & \Omega_{2}(5) = \lcb z \in \complex \, : \, |z^{3} + 2 z^{2} - (1-i) z - 11| \leq 2.2720 \rcb \; .
\end{eqnarray*}
Figure~\ref{Pelletfig} shows $\Omega_{1}(5)$ in light gray and $\Omega_{2}(5)$ in dark gray for $q_{A}$, while
its zeros are represented by the white dots. The two circles are the boundaries of the disks $O(0;r)$ and $O(0;R)$. 
The smaller
disk, which contains five zeros, is contained in $\Omega_{1}(5)$, while $\Omega_{2}(5)$, which contains the
remaining three zeros, lies outside the larger disk.                  
The disk $\obar (0;r) $ is clearly preferable to $\Omega_{1}(5)$,
whereas $\Omega_{2}(5)$ is clearly preferable to  $O^{c}(0;R)$, as predicted by the theorem.
We remark that $\Omega_{2}(k)$ does not necessarily consist of disjoint regions, although it often happens,
as in this case. The boundary of $\Omega_{2}(5)$ has foci at $1.8291 - 0.1119i$, $-2.1035 + 1.5937i$, 
and $-1.7257 - 1.4818i$, which can be enclosed in disjoint disks of radius $1.5$ centered at these foci. 
Since $(1.5)^3=3.3750$ and $2.2720 < 3.3750$, $\Omega_{2}(5)$ must consist of three disjoint regions with 
simple boundaries, as explained at the end of Section 2.
By part (3) of Theorem~\ref{geometricPellet} each must contain exactly one zero of $q_{A}$.

%
%
%
%

\begin{figure}[H]
\begin{center}
\includegraphics[width=0.45\linewidth]{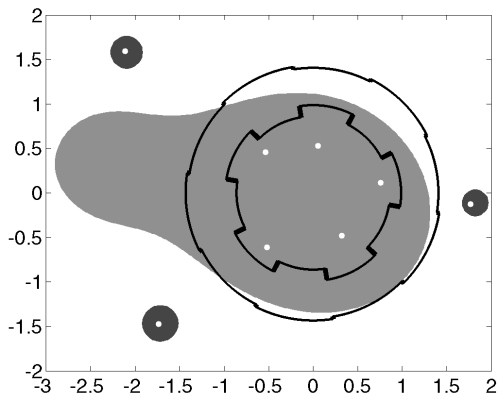}
\caption{The sets $\Omega_{1}(5)$ (light gray) and $\Omega_{2}(5)$ (dark gray) for $q_{A}$, compared to 
the refined Pellet region.}
\label{refinedPelletfig}
\end{center}
\end{figure}

In~\cite{MardenPaper} and~\cite[Theorem (29,1), p.130 and Exercise 2, p. 133]{Marden} a refinement of
Pellet's theorem was derived, requiring the solution of an additional equation, related to 
$f_{k}(x)=0$ in Theorem~\ref{geometricPellet}, along with a nonzero requirement on one additional coefficient. 
It was further generalized 
in \cite[Theorem 1]{Levinger}. The refinement
leads to a slightly better gear-wheel shaped region (instead of an annulus) that does not contain any zeros of the polynomial.
In Figure~\ref{refinedPelletfig} we have compared $\Omega_{1}(5)$ and $\Omega_{2}(5)$ for $q_{A}$
to the refined region from~\cite[Theorem (29,1), p.130]{Marden}.
The radii of the circles determining the inner boundary of this region 
are $0.8612$ and $0.9872$, whereas for the outer boundary they are given by $1.4065$ and $1.4331$.
%
%
%
%

\begin{figure}[H]
\begin{center}
\includegraphics[width=0.40\linewidth]{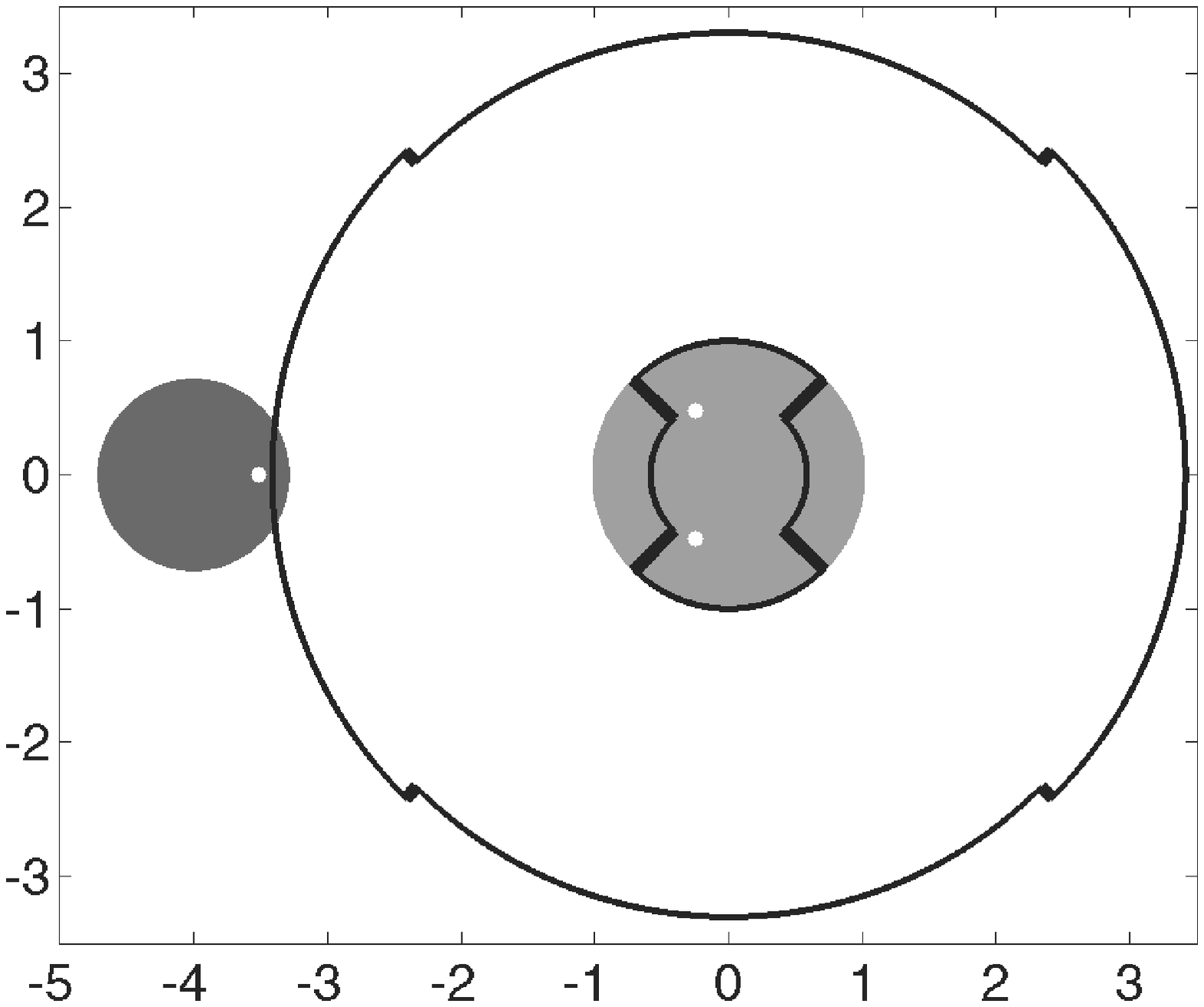}
\;\;\;\;
\;\;\;\;
\includegraphics[width=0.40\linewidth]{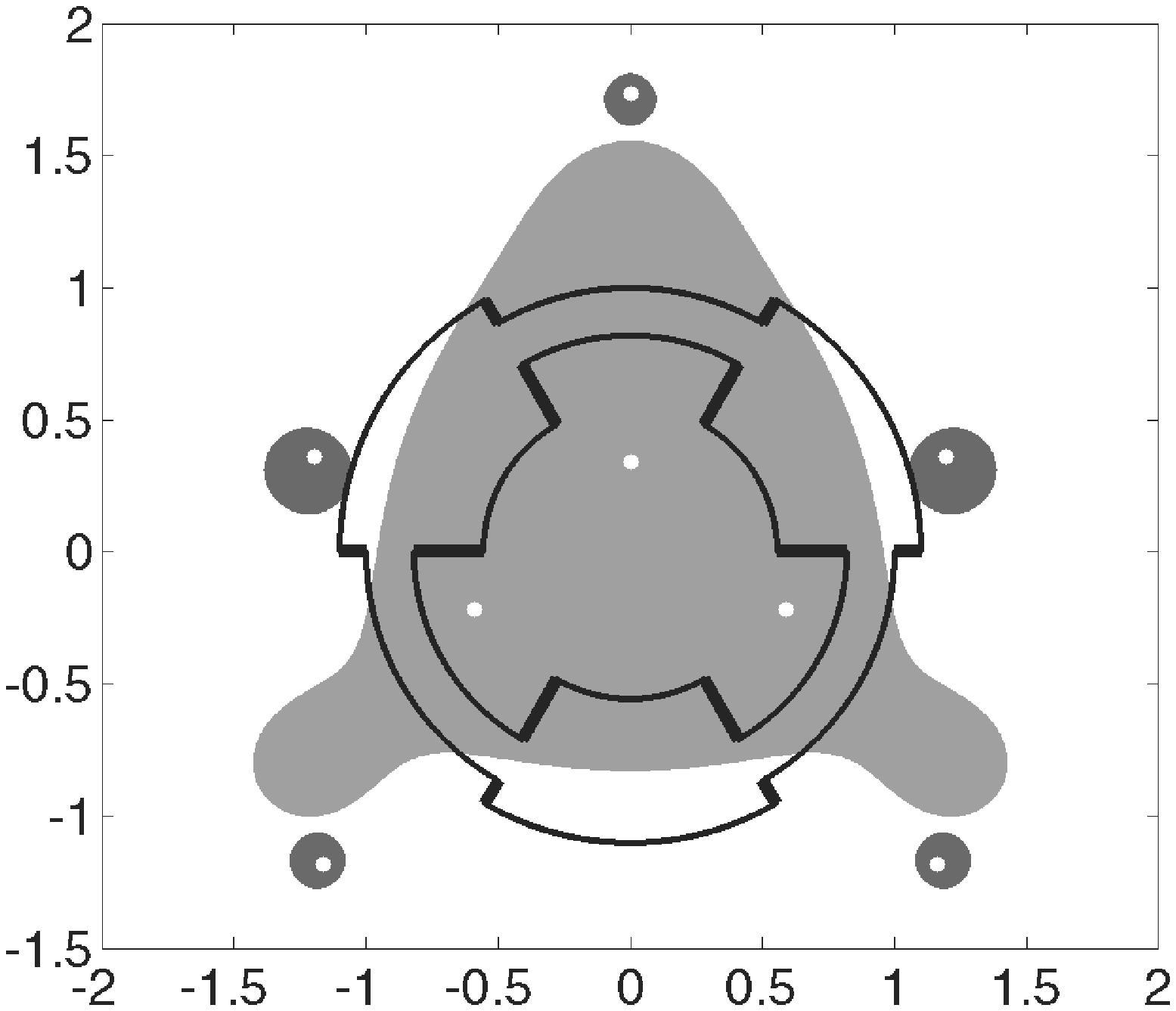}
\caption{The sets $\Omega_{1}(k)$ (light gray) and $\Omega_{2}(k)$ (dark gray) for $q_{B}$ and $q_{C}$ with $k=2$
and $k=3$, respectively, compared to the refined Pellet region.}
\label{q2q3fig}
\end{center}
\vskip -0.5cm
\end{figure}

\begin{figure}[H]
\begin{center}
\includegraphics[width=0.45\linewidth]{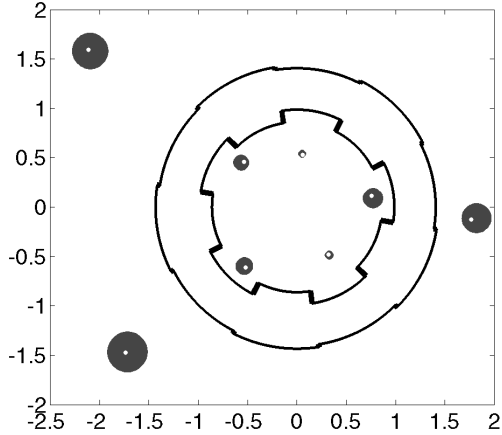}
\caption{Inclusion sets obtained from both the polynomial and its reciprocal for $q_{A}$.}
\label{recipfig}   
\end{center}
\vskip -0.5cm
\end{figure}

Let us consider two more examples, namely, $q_{B}=z^3+4z^2+2z+1$ (the example in~\cite{Levinger}) for $k=2$, 
and $q_{C}=z^{8}+3iz^5+z^4-8iz^3+2iz+1$ for $k=3$. Similar conclusions as for $q_{A}$ can be drawn for these polynomials
as can be seen from  
Figure~\ref{q2q3fig}, which shows (using the same conventions as before)
the corresponding sets $\Omega_{1}(k)$ and $\Omega_{2}(k)$ for $q_{B}$ (left) and 
$q_{C}$ (right) with $k=2$ and $k=3$, respectively, together with the refined Pellet region from 
~\cite[Theorem (29,1), p.130]{Marden}. 
We leave out the details for brevity.
%
%
%
%
$ $ \newline {\bf Remarks.}
{\bf (1)} By using the reciprocal polynomial, similar lemniscates can be derived for the reciprocals
of the zeros of a polynomial, leading to corresponding inclusion regions for the zeros themselves,
although these are more complicated. They can be combined with our previous results
as is illustrated for $q_{A}$ in Figure~\ref{recipfig},
which also shows the same refined Pellet region as before. Results of a very similar nature are obtained 
for $q_{B}$ and $q_{C}$.
%
%
%
%
%
$ $ \newline {\bf (2)} A converse of Pellet's theorem is stated in~\cite{Walsh} 
(see also \cite[Theorem (28,3), p.129]{Marden}), whose proof was later corrected in~\cite{Ostrowski}.
It can be formulated as follows.
Let $a_{j}$ ($j=1,2,...,n$) be fixed complex coefficients and $\epsilon_{j}$ ($j=1,2,...,n$) be arbitrary 
complex numbers with $|\epsilon_{j}|=1$ for all $j=1,2,...,n$.
Let $\tau$ be any positive number such that it is not a zero of any polynomial 
$\sum_{j=1}^{n} a_{j} \epsilon_{j} z^{j}$
and let every polynomial of that form have $k$ zeros ($0 < k < n$) in $\obar(0;\tau)$.
Then $f_{k}$ (with $f_{k}$ as defined in Theorem~\ref{Pellet}) has two positive roots $r$ and $R$
such that $r < \rho < R$.
%

It may be possible to construct a similar converse of Theorem~\ref{geometricPellet}, 
although its proof would be well beyond the scope of this paper.
$ $ \newline {\bf (3)} A sufficient condition for the $n-k$ largest zeros to have a modulus strictly larger 
than $R$ is for $\Omega_{2}(k)$ to be contained in the interior of $O^{c}(0;R)$. A similar condition based
on the reciprocal polynomial holds for the $k$ smallest zeros. In~\cite{ZhangXuNiu}, different (but more explicit)  
conditions were derived for the same situation. The special case $r=R=\rho$ was considered in~\cite{Walsh} (also 
mentioned in~\cite[Theorem (28,2), p.129]{Marden}), where it was shown that the polynomial then has $\ell$ 
($\ell \geq 0$) double roots on the circle $|z|=\rho$, $k-\ell$ zeros inside and $n-k-\ell$ zeros outside this
circle, respectively. Theorem~\ref{geometricPellet} adds to this result a sufficient condition 
guaranteeing that no double roots lie on the circle $|z|=\rho$, as follows:
from its proof, we have in this case that the two
disks in the Gershgorin set of $D_{x}^{-1}M_{n-k}(p)D_{x}$ are tangent to each other when 
$p_{n-k}(\zeta) = - P_{n-k}(\rho) \ak/|\ak|$ for a point $\zeta \in \Omega_{1}(k) \cap \Omega_{2}(k)$, 
i.e., when 
\beq 
\label{intsecpoints}
p_{n-k}(\zeta) + \dfrac{P_{n-k}(\rho)}{|\ak|} \, \ak = 0 \; .
\eeq   
This means that $\Omega_{1}(k)$ and $\Omega_{2}(k)$ touch at exactly $n-k$ points,
namely the zeros of the polynomial of degree $n-k$ in~(\ref{intsecpoints}).
If those points do not on the circle $|z|=\rho$, then, 
because the zeros of a polynomial are continuous functions of their coefficients, none in the 
group of largest zeros of $p$ can cross the gap between $\Omega_{2}(k)$ and the circle when the coefficients are 
continuously perturbed 
from a situation of two very close distinct values for $r$ and $R$ to one where $r=R$.
In view of~\cite{Walsh}, this means that in such a case no double roots of $p$ can lie on the circle.


$ $ \newline {\bf (4)} The set $\Omega_{2}(k)$ can also be derived without using Gershgorin's theorem, in the following
way. The equation $p(z)=0$ can be written as 
\bdis
z^{k} \lb z^{n-k} + a_{n-1}z^{n-k-1} + \dots + a_{k} \rb = -a_{k-1}z^{k-1} - \dots -a_{0} \; ,
\edis
from which we have
\bdis
\left | z^{n-k} + a_{n-1}z^{n-k-1} + \dots + a_{k} \right | 
\leq \dfrac{|a_{k-1}|}{|z|} + \dots + \dfrac{|a_{0}|}{|z|^{k}} \; \cdot
\edis
The zeros $\{ z_{i} \}_{i=k+1}^{n}$, which satisfy $|z_{i}| \geq \rtwo$, must then also satisfy
\bdis
\left | z_{i}^{n-k} + a_{n-1}z_{i}^{n-k-1} + \dots + a_{k} \right | 
\leq \dfrac{|a_{k-1}|}{\rtwo} + \dots + \dfrac{|a_{0}|}{\rtwo^{k}} \; ,          
\edis
and therefore $z_{i} \in \Omega_{2}(k)$ for $k+1 \leq i \leq n$.

On the other hand, Gershgorin's theorem permits a unified and convenient treatment of both the previous and 
the next theorem, our second result, which presents a new generalization of Theorem~\ref{Cauchy}.

\begin{theorem}
\label{genCauchy}
Let $p(z)=z^{n} + a_{n-1} z^{n-1} + \dots + \aone z + \azero$ be a polynomial with complex coefficients, 
$\azero \neq 0$, $1 \leq k \leq n-1$, $n \geq 3$, and with zeros $\{z_{i}\}_{i=1}^{n}$.
Let $\{p_{j}\}_{j=1}^{n-1}$ be the associated polynomials of $p$, 
let $P_{j}$ be the polynomial obtained from $p_{j}$ by replacing its coefficients with their moduli, and let
\bdis
\mu(k,x)=\sum_{j=0}^{k-1} |a_{j}|x^{j-k} \; .  
\edis
For $j=0,1,\dots,n-1$, let $s_{j}$ be the unique positive root of 
\bdis
h_{j}(x)=x^{n} + |a_{n-1}| x^{n-1} + \dots + |a_{j+1}| x^{j+1} - |a_{j}| x^{j} - |a_{j-1}|x^{j-1} - \dots  - |\azero| \; .
\edis

Then 
\newline {\bf (1)} all the zeros of $p$ are contained in $\Upsilon_{1}(k)$ and also in $\Upsilon_{2}(k)$,
where
\begin{eqnarray*}
& & \Upsilon_{1}(k) = \lcb z \in \complex \, : \, |p_{n-k}(z) | \leq P_{n-k}(s_{k}) = \mu(k,s_{k}) + |\ak| \rcb  \; , \\    
& & \Upsilon_{2}(k) = \lcb z \in \complex \, : \, |p_{n-k}(z) + a_{k}| \leq \mu(k,s_{k-1}) = P_{n-k}(s_{k-1}) + |\ak| \rcb \; ,
\end{eqnarray*}
each of which is the closed interior of a lemniscate;
\newline {\bf (2)} if $\Upsilon_{1}(k)$ or $\Upsilon_{2}(k)$ consists of disjoint regions whose boundaries are 
simple closed (Jordan) curves and $\ell$ is the number of foci of the corresponding lemniscate contained in any 
such region, then that region contains $\ell$ zeros of $p$ when that region does not contain the origin, and when 
it does contain the origin, then it contains $\ell + k$ zeros of $p$.
\end{theorem}
\prf
We begin by observing that $h_{j}$ can be written as
\beq
\label{hmueq}
h_{j}(x) = x^{j} \lb P_{n-j}(x) -|a_{j}| - \mu(j,x) \rb \; .
\eeq
Let us now consider once more the matrix
$D^{-1}_{\ix} M_{n-k}(p) D_{\ix}$ for any $x > 0$, where $D_{\ix}$ is as in Lemma~\ref{lemmaMksimilarity},
whose eigenvalues are the numbers $\{p(z_{i})\}_{i=1}^{k}$. Its Gershgorin column set, which contains these
numbers, is the union of the two disks
$\obar \lb 0;P_{n-k}(x) \rb $ and $\obar \lb -a_{k};\mu(k,x) \rb $. As $x$ increases, so does $P_{n-k}(x)$,
while $\mu(k,x)$ decreases. The former disk will then encompass the latter when  
$P_{n-k}(x) = |a_{k}| + \mu(k,x)$, or, as can be seen from~(\ref{hmueq}), $h_{k}(x)=0$. That is, when $x=s_{k}$. 
All the zeros of $p$ are then contained in the set
\bdis
\Upsilon_{1}(k) = \lcb z \in \complex \, : \, |p_{n-k}(z) | \leq P_{n-k}(s_{k}) \rcb  \; .
\edis
Because $h_{k}(s_{k})=0$, the right-hand side of the inequality defining $\Upsilon_{1}(k)$ can be replaced 
by~$\mu(k,s_{k})~+~|\ak|$.

On the other hand, we can let $x$ decrease until the disk centered at $-a_{k}$ encompasses the one centered 
at the origin. This happens when 
$\mu(k,x) = |a_{k}| + P_{n-k}(x)$. Since, by using~(\ref{hmueq}), we can write
\begin{eqnarray*}
P_{n-k}(x) + |a_{k}|  - \mu(k,x) & = & x^{-1} \Bigl( x \bigl (  P_{n-k}(x) + |a_{k}| \bigr ) - x\mu(k,x) \Bigr ) \\  
                                    & = & x^{-1} \Bigl ( P_{n-k+1}(x) - \bigl ( |a_{k-1}| + \mu(k-1,x) \bigr ) \Bigr ) \\
                                    & = & x^{-k} h_{k-1}(x) \; , \\
\end{eqnarray*}
$ $ \newline \vskip -1.50cm
\hskip -0.5cm we conclude that $x=s_{k-1}$. All the zeros of $p$ are then contained in the set
\bdis
\Upsilon_{2}(k) = \lcb z \in \complex \, : \, |p_{n-k}(z) + a_{k} | \leq \mu(k,s_{k-1}) \rcb  \; .
\edis
Because $h_{k-1}(s_{k-1})=0$, the right-hand side of the inequality defining $\Upsilon_{1}(k)$ can be replaced 
by~$P_{n-k}(s_{k-1}) + |\ak|$. 

If there exists a disjoint region $\scr$ of $\Upsilon_{1}(k)$ with a simple boundary, we proceed
similarly as in the proof of Theorem~\ref{geometricPellet} and first define 
\bdis
\widetilde{\Upsilon}_{1}(k) = \lcb z \in \complex \, : \, |p_{n-k}(z) | \leq P_{n-k}(s)  \rcb  \; ,   
s > s_{k} \; ,
\edis
where $s-s_{k}$ is small enough so that $\widetilde{\Upsilon}_{1}(k)$ (which contains $\Upsilon_{1}(k)$
because $P_{n-k}(s)~>~P_{n-k}(s_{k})$) 
has the same number of 
disjoint regions with simple boundary as $\Upsilon_{1}(k)$, one of which, $\widetilde{\scr}$, contains $\scr$. 
Any $z \in \partial \widetilde{\scr}$ satisfies 
\bdis
P_{n-k}(|z|) \geq |p_{n-k}(z)| = P_{n-k}(s) \; ,
\edis
implying that $|z| \geq s$ because $P_{n-k}$ is increasing for positive arguments.
In addition, because $s > s_{k}$, we also have that $P_{n-k}(s) > \mu(k,s) + |\ak|$. Then we can write 
for any $z \in \partial \widetilde{\scr}$:
\bdis
|p_{n-k}(z)| = P_{n-k}(s) > |\ak| + \dfrac{|\akmo|}{s} + \dots + \dfrac{|\azero|}{s^{k}}  
             \geq |\ak| + \dfrac{|\akmo|}{|z|} + \dots + \dfrac{|\azero|}{|z|^{k}} \; ,
\edis
from which it follows that 
\bdis
|z|^{k}|p_{n-k}(z)| >  |\ak||z|^{k} + |\akmo||z|^{k-1} + \dots + |\azero| \; ,  
\edis
and therefore
\bdis |z^{k}p_{n-k}(z)| > \left | \ak z^{k} + \akmo z^{k-1} + \dots + \azero \right | \; .  
\edis
Since we assumed that $\partial \widetilde{\scr}$ is a simple closed (Jordan) curve, 
we conclude from Rouch\'{e}'s theorem 
that the polynomial $p$ 
has as many zeros in $\widetilde{\scr}$ as the polynomial $z^{k}p_{n-k}(z)$. If $0 \notin \widetilde{\scr}$, 
then that number is the 
number of foci of the lemniscate forming the boundary of $\widetilde{\Upsilon}_{1}(k)$ that lie in $\widetilde{\scr}$, 
i.e., those zeros of $p_{n-k}$ that are contained in $\widetilde{\scr}$. 
If $0 \in \widetilde{\scr}$, then $k$ is added to the number of foci. 
Because this remains true as $s \rightarrow s_{k}^{+}$, the proof follows. 
The proof for $\Upsilon_{2}(k)$ is analogous. \qed

Following are the four lowest-order lemniscates containing all the zeros of $p$:
\begin{eqnarray*}
& & \Upsilon_{1}(n-1) = \lcb z \in \complex \, : \, |z| \leq s_{n-1} \rcb \; , \\ 
& & \Upsilon_{2}(n-1) = \lcb z \in \complex \, : \, |z + a_{n-1}| \leq |a_{n-1}| + s_{n-2} \rcb \; , \\
& & \Upsilon_{1}(n-2) = \lcb z \in \complex \, : \, |z(z+a_{n-1})| \leq s_{n-2}(s_{n-2}+|a_{n-1}|) \rcb \; , \\ 
& & \Upsilon_{2}(n-2) = \lcb z \in \complex \, : \, |z(z + a_{n-1}) + a_{n-2}| 
    \leq s_{n-3}(s_{n-3}+|a_{n-1}|) + |a_{n-2}| \rcb \; . \\
\end{eqnarray*}

\vskip -0.5cm
The sets $\Upsilon_{1}(n-1)$ and $\Upsilon_{2}(n-1)$ are closed disks, whereas the sets
$\Upsilon_{1}(n-2)$ and $\Upsilon_{2}(n-2)$ are the closed interiors of ovals of Cassini 
(\cite[p.153-155]{Lawrence}). 
The zeros of $p$ lie in the intersection of all of the aforementioned
sets, although, with the exception of the sets that are disks, such an intersection is in general 
difficult to compute.
When all the coefficients of the polynomial are real, then all our 
inclusion sets are symmetric with respect to the real axis.

\begin{figure}[H]
\begin{center}
\hskip -0.5cm 
\includegraphics[width=0.45\linewidth]{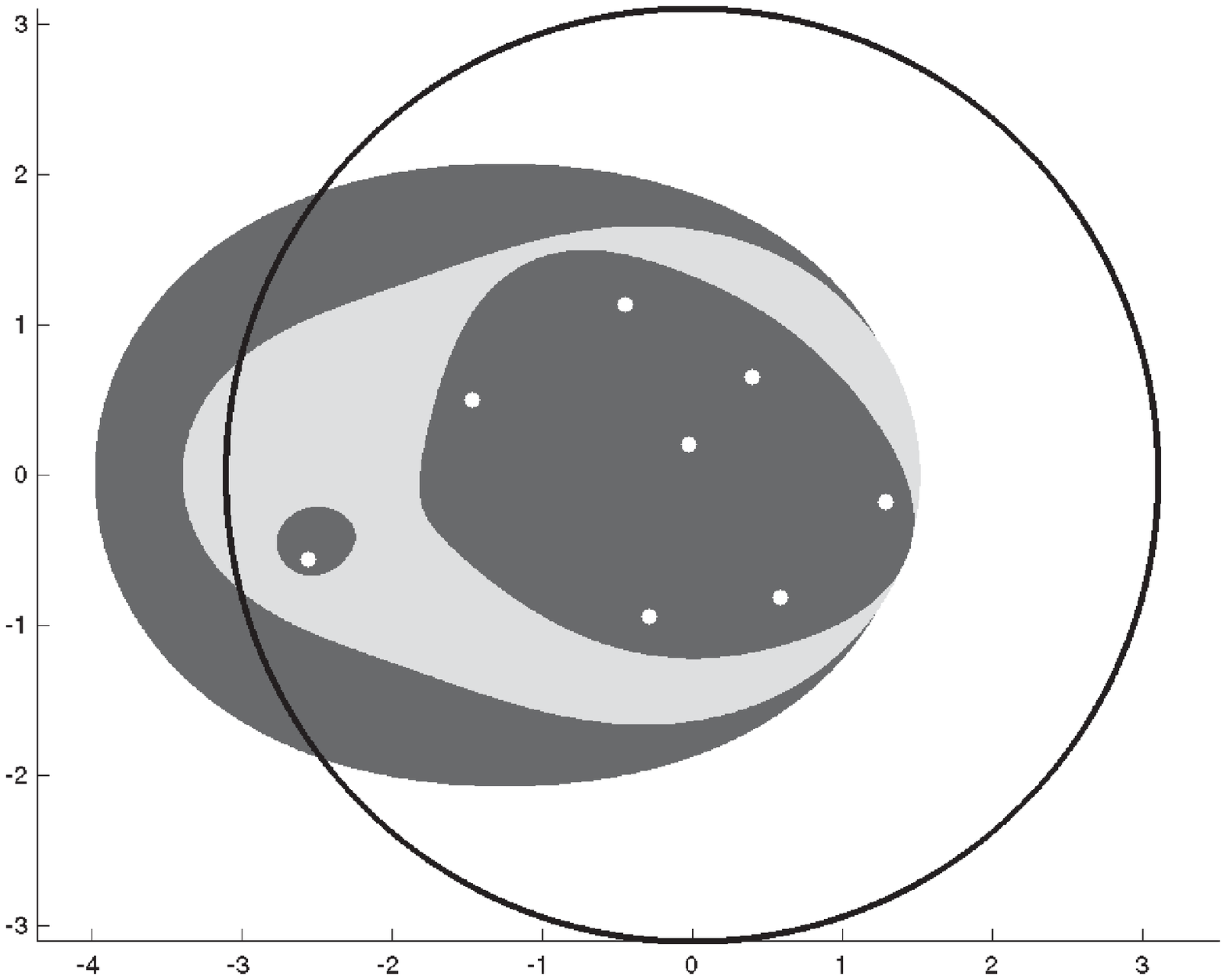}
\hskip 1.5cm   
\includegraphics[width=0.45\linewidth]{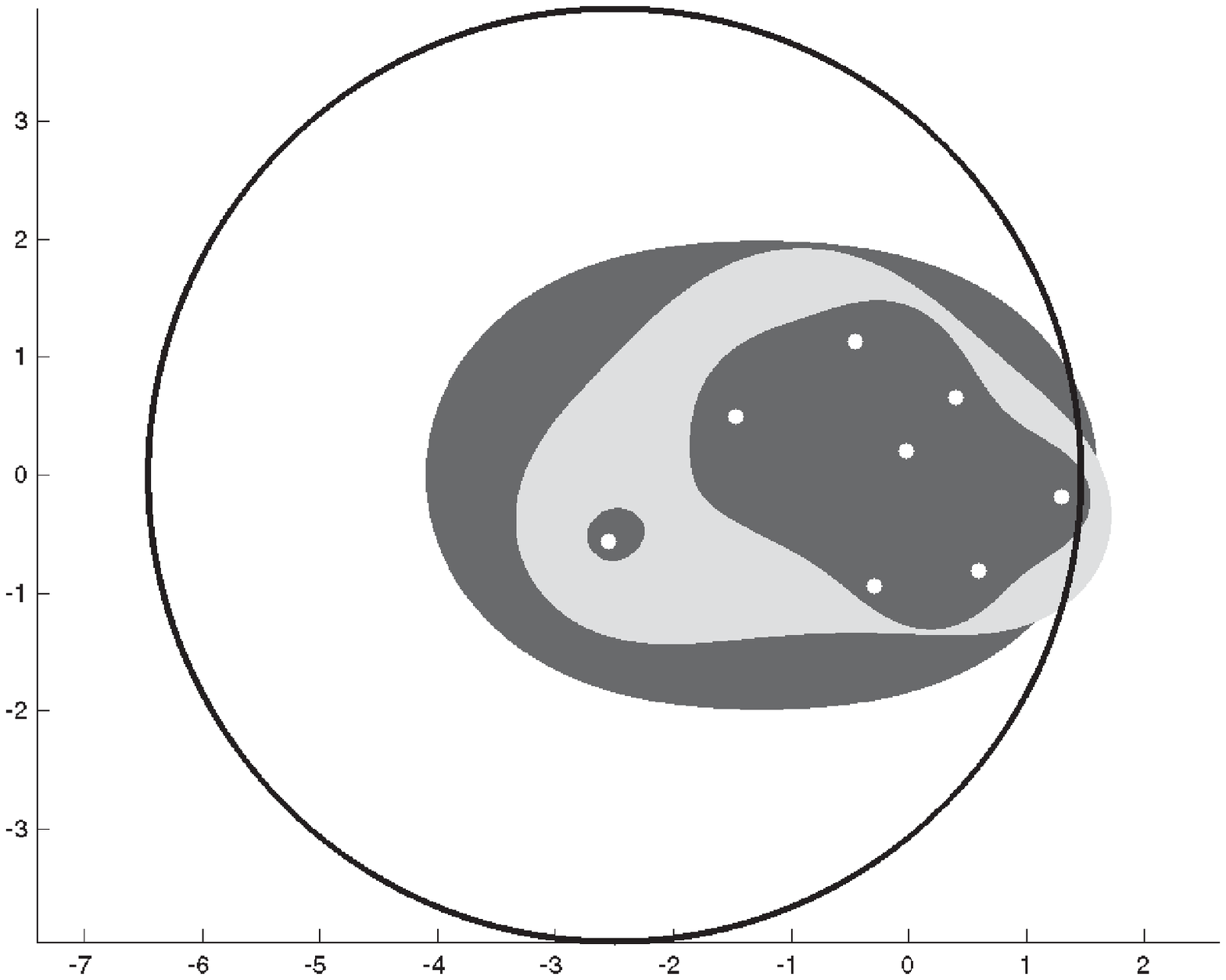}
\caption{The sets $\Upsilon_{1}(k)$ (left) and $\Upsilon_{2}(k)$ (right) for $q_{D}$ with $k=6,5,3$.}
\label{Cauchyfig}   
\end{center}
\end{figure}

To illustrate Theorem~\ref{genCauchy}, we consider the polynomial $q_{D}$, defined by
\bdis
q_{D}(z) = 
z^{8} + \dfrac{5}{2} z^{7} - \dfrac{1}{2} z^{6} - (2-3i) z^{5} - \dfrac{1}{2} z^{4} - 4 z^{3} + 2 z^{2} - 5 z + i .
\edis
Figure~\ref{Cauchyfig} shows the corresponding sets $\Upsilon_{1}(k)$ on the left and $\Upsilon_{2}(k)$ on the 
right for $k=6,5,3$. To better compare them, the lemniscates have been superimposed in alternating shades of light 
and dark gray.  
The circles mark the boundaries of the sets $\Upsilon_{1}(n-1)=\Upsilon_{1}(7)$ (left) and 
$\Upsilon_{2}(n-1)=\Upsilon_{2}(7)$ (right).
The largest proper lemniscates in dark gray, the ones in light gray,
and the smallest ones in dark gray correspond to $k=6$, $k=5$, and $k=3$, respectively.
The white dots are the zeros of $q_{D}$.
%
%
%
%

$ $ \newline {\bf Remarks.}
{\bf (1)} 
The sets $Y_{1}(k)$ and $Y_{2}(k)$ are not necessarily nested for successive values
of $k$, although they do tend to become smaller as $k$ decreases.
\newline {\bf (2)}
Theorem~\ref{genCauchy} is a generalization of Theorem~\ref{Cauchy},
since that theorem is obtained for $\Upsilon_{1}(k)$ with $k=n-1$ in Theorem~\ref{genCauchy}.
Lemniscates are more interesting geometric regions than disks, 
and even low-order lemniscates
can already provide significantly improved inclusion regions for the zeros of a polynomial. 
The computation of the positive root of any $h_{j}$ is inexpensive compared to the computation of 
the (generally complex) zeros of $p$.
\newline {\bf (3)}
The last part of the previous theorem's proof shows that it can be proven 
with just Rouch\'{e}'s theorem, without using Gershgorin's theorem. However, the latter provides
a natural explanation for the appearance of the functions $h_{j}$ that would otherwise be 
lacking, and also creates a larger framework in which both Pellet's and the generalized Cauchy 
theorems are obtained as special cases for special values of $x$ in $D_{x}^{-1} M_{n-k}(p)D_{x}$.
This framework generates infinitely many other inclusion regions depending on the values of $x$, e.g.,
the value for which both disks in the proof of Theorem~\ref{genCauchy} have the same radius, leading to
a union of the interiors of two lemniscates, to give but one example.
\newline {\bf (4)}
The $\Upsilon_{1}(k)$ and $\Upsilon_{2}(k)$ sets can sometimes be simplified if the coefficients of the polynomial 
exhibit certain patterns. For instance, if the leading
coefficients fit the pattern of a power of $(z-a)$ for some complex number a, then the corresponding
set, defined by a polynomial of the same degree as the power, becomes a simple disk. Consider
as an example the set $\Upsilon_{2}(3)$ for the polynomial $z^5 + 2iz^4 - z^3 + z^2 + 3z - 1$ , whose first three
coefficients are the same as those of $(z + i)^2$. It is given by
\bdis
\Upsilon_{2}(3) = \lcb z \in \complex : |z + i| = s_{2} + 1 \rcb \; .
\edis
The inclusion sets for lacunary polynomials with several consecutive leading zero coefficients can be treated in 
the same way since their leading coefficients fit the pattern of a power of $(z-0)$.


\end{document}